\documentclass[12pt,reqno]{article}
\usepackage{graphicx}
\usepackage{amsmath}
\usepackage{amsthm}
\usepackage{latexsym,bm}
\usepackage{amssymb}
\usepackage{indentfirst}
\newtheorem{thm}{Theorem}[section]

\newtheorem{lem}[thm]{Lemma}
\newtheorem{prop}[thm]{Proposition}

\numberwithin{equation}{section}
\begin{document}

\title{{\bf On complete hypersufaces with constant mean and scalar curvatures in Euclidean spaces}}

\author{{\bf Roberto Alonso N\'u\~nez}}

\date{}
\maketitle
\footnotetext[1]{2010 {\it Mathematics Subject Classication.} Primary: 53C40; Secondary: 53C42.}

\begin{quote}
\small {\bf Abstract}.
Generalizing a theorem of Huang, Cheng and Wan classified the complete hypersurfaces of $\mathbb R^4$ with non-zero constant mean curvature and constant scalar curvature. In our work, we obtain results of this nature in higher dimensions. In particular, we prove that if a complete hypersurface of $\mathbb R^5$ has constant mean curvature $H\neq 0$ and constant scalar curvature $R\geq\frac{2}{3}H^2$, then $R=H^2$, $R=\frac{8}{9}H^2$ or $R=\frac{2}{3}H^2$. Moreover, we characterize the hypersurface in the cases $R=H^2$ and $R=\frac{8}{9}H^2$, and provide an example in the case $R=\frac{2}{3}H^2$. The proofs are based on the principal curvature theorem of Smyth-Xavier and a well known formula for the Laplacian of the squared norm of the second fundamental form of a hypersurface in a space form.
\end{quote}

\section{Introduction}
\vskip10pt

A well known result of Klotz and Osserman \cite{Klot} states that the round spheres and the circular cylinders are the only complete surfaces in the Euclidean $3$-space $\mathbb R^3$ with non-zero constant mean curvature whose Gaussian curvature does not change sign. In higher dimensions, Nomizu and Smyth \cite{Nom_Smyth} proved the following result: {\it If $M^n$ is a complete Riemannian manifold with non-negative sectional curvature and constant scalar curvature, and $f:M^n\to\mathbb R^{n+1}$ is an isometric immersion with constant mean curvature, then $f(M^n)$ is a generalized cylinder $\mathbb R^{n-k}\times\mathbb{S}^{k}$, for some $k=1,\ldots,n$.}

Cheng and Yau proved in \cite{ChengYauMean} that the above mentioned result of Nomizu and Smyth holds without the assumption that the scalar curvature is constant, and in \cite{ChengYauScal} that it holds without the assumption that the mean curvature is constant. Extending these results of Cheng and Yau, Hartman proved in \cite{HART} that if $M^n$ is a complete Riemannian manifold with non-negative sectional curvature and $f:M^n\to\mathbb R^{n+1}$ is an isometric immersion with positive constant $r$-th mean curvature $H_r$, for some $r=1,...,n$, then $f(M^n)=\mathbb R^{n-d}\times\mathbb S^{d}$, for some $r\leq d\leq n$.

In view of the above discussion, a question that arises naturally is whether the theorem of Nomizu-Smyth mentioned in the first paragraph holds without the assumption that the sectional curvatures of $M^n$ are non-negative. For n = 3, a partial answer to this question was given by Huang \cite{HUANG}:

\vskip7pt

\noindent{\bf Teorema} (Huang).
{\em Let $M^3$ be a complete and connected Riemannian manifold with non-negative constant scalar curvature $R$, and let $f:M^3\to\mathbb R^4$ be an isometric immersion with non-zero constant mean curvature $H$. Then, $R=H^2$, $R=\frac{3}{4}H^2$ or $R=0$. When $R=H^2$, $f(M^3)=\mathbb S^3\big(\frac{1}{|H|}\big)$, and when $R=\frac{3}{4}H^2$, $f(M^3)=\mathbb R\times \mathbb S^2\big(\frac{2}{3|H|}\big)$.}

\vskip7pt

Later, Cheng and Wan \cite{ChengWan} improved the above theorem of Huang in two aspects. They not only showed that the hypothesis $R\geq0$ is superfluous, but also characterized the hypersurface in the case $R=0$. More precisely, they obtained the following result.

\vskip7pt

\noindent{\bf Teorema} (Cheng-Wan).
{\em Let $M^3$ be a complete and connected Riemannian manifold with constant scalar curvature $R$, and let $f:M^3\to\mathbb R^4$ be an isometric immersion with non-zero constant mean curvature $H$. Then $R=H^2$, $R=\frac{3}{4}H^2$ or $R=0$. When $R=H^2$, $f(M^3)=\mathbb S^3\big(\frac{1}{|H|}\big)$; when $R=\frac{3}{4}H^2$, $f(M^3)=\mathbb R\times \mathbb S^2\big(\frac{2}{3|H|}\big)$ and when $R=0$, $f(M^3)=\mathbb R^2\times \mathbb S^1\big(\frac{1}{3|H|}\big)$.}

\vskip7pt

In the same spirit of the results of Huang and Cheng-Wan referred to above, we establish in this work the following theorems:
\begin{thm}\label{teorema1}
Let $M^4$ be a complete and connected Riemannian manifold with constant scalar curvature $R$, and $f:M^4\to\mathbb R^5$ an isometric immersion with non-zero constant mean curvature $H$. If $R\geq\frac{2}{3}H^2$ then $R=H^2$, $R=\frac{8}{9}H^2$ or $R=\frac{2}{3}H^2$. When $R=H^2$, $f(M^4)=\mathbb S^4\big(\frac{1}{|H|}\big)$ and when $R=\frac{8}{9}H^2$,
$f(M^4)=\mathbb R\times\mathbb S^3\big(\frac{3}{4|H|}\big)$.
\end{thm}
\noindent{\bf Remark 1.} When $R=\frac{2}{3}H^2$, one has as an example the cylinder $\mathbb R^2\times \mathbb S^2\big(\frac{1}{2|H|}\big)$.

\begin{thm}\label{teorema2}
Let $M^5$ be a complete and connected Riemannian manifold with constant scalar curvature $R$, and $f:M^5\to\mathbb R^6$ an isometric immersion with non-zero constant mean curvature $H$ and non-negative $4$-th mean curvature $H_4$. If $R\geq\frac{5}{8}H^2$ then $R=H^2$, $R=\frac{15}{16}H^2$, $R=\frac{5}{6}H^2$ or $R=\frac{5}{8}H^2$. When $R=H^2$, $f(M^5)=\mathbb S^5\big(\frac{1}{|H|}\big)$ and when $R=\frac{15}{16}H^2$, $f(M^5)=\mathbb R\times\mathbb S^4\big(\frac{4}{5|H|}\big)$.
\end{thm}
\noindent{\bf Remark 2.} When $R=\frac{5}{6}H^2$, one has as an example the cylinder $\mathbb R^2\times\mathbb S^3(\frac{3}{5|H|})$ and when $R=\frac{5}{8}H^2$, the cylinder $\mathbb R^3\times\mathbb S^2(\frac{2}{5|H|})$.
\begin{thm}\label{teorema3}
Let $M^n$ be a complete and connected Riemannian manifold of dimension $n\geq 3$ and constant scalar curvature $R$, and let $f:M^n\to\mathbb R^{n+1}$ be an isometric immersion with non-zero constant mean curvature $H$. If $HH_3\geq0$ and $\;0\leq R\leq\frac{nH^2}{2(n-1)}$, then $R=0$ or $R=\frac{nH^2}{2(n-1)}$.
In case that $R=0$, $f(M^n)=\mathbb R^{n-1}\times\mathbb S^1\big(\frac{1}{n|H|}\big)$ and in case that $R=\frac{nH^2}{2(n-1)}$, $f(M^n)=\mathbb R^{n-2}\times\mathbb S^2\big(\frac{2}{n|H|}\big)$.
\end{thm}
\noindent{\bf Remark 3.} If it were possible to remove the hypothesis $HH_3\geq0$ in Theorem \ref{teorema3}, then a combination of this theorem (for $n=4$) with Theorem \ref{teorema1} would provide an extension of the theorem of Cheng-Wan for hypersurfaces with non-negative scalar curvatures in $\mathbb R^5$.

\vskip10pt

\noindent{\bf Acknowledgements.} This work is part of the author's doctoral thesis at Universidade Federal Fluminense. The author would like to express his sincere gratitude to his advisor Prof. Francisco Fontenele for the continuous support of his Ph.D study and research, and to CAPES (Brazil) for the financial support.

\vskip10pt
\section{Preliminaries}\label{preliminar}
\vskip10pt

Given an isometric immersion $f:M^n\rightarrow N_c^{n+1}$ of an orientable $n$-dimensional Riemannian manifold $M^n$ into an orientable $(n+1)$-dimensional Riemannian manifold $N^{n+1}_c$ of constant sectional curvature $c$, we denote by $A$ the shape operator of $f$ with respect to a global unit normal vector field $\xi$, and by $\lambda_1,\ldots,\lambda_n$ the eigenvalues of $A$ (the principal curvatures of $M^n$). It is well known that if we label the principal curvatures at each point by the condition $\lambda_1\leq\ldots\leq\lambda_n$, then the principal curvature functions $\lambda_i:M\rightarrow\mathbb R$, $i=1,\ldots,n$, become  continuous.

The $r$-th mean curvature $H_r$, $1\leq r\leq n$, of the immersion is defined by
\begin{eqnarray}\label{Hr}
{n \choose r}H_r=S_r:=\sum_{i_1<\ldots<i_r}\lambda_{i_1}\ldots \lambda_{i_r}.
\end{eqnarray}
Notice that $H_1$ is the mean curvature $H$ and
$H_n=\lambda_1\lambda_2\ldots\lambda_n$ is the Gauss-Kronecker curvature of the immersion. The function $H_2$ is up to a constant the (normalized) scalar curvature $R$ of $M^n$. Indeed, if for a given $p\in M$ we consider an orthonormal basis $\{e_1,\ldots,e_n\}$ of $T_pM$ that diagonalizes $A$, then the sectional curvature $K(e_i,e_j)$ of the plane spanned by $e_i$ and $e_j$ is, by the Gauss equation, given by
\begin{eqnarray}\label{tb117}
K(e_i,e_j)=c+\lambda_i\lambda_j,
\end{eqnarray}
and so
\begin{eqnarray}\label{hr2}
R=\frac{1}{{n \choose 2}}\sum_{i<j}K(e_i,e_j)=\frac{1}{{n \choose 2}}\sum_{i<j}(c+\lambda_i\lambda_j)=c+H_2.
\end{eqnarray}

For what follows it is convenient to write the $r$-th mean curvature $H_r(p)$ of $M^n$ at a point $p$ as
\begin{eqnarray}\label{hr3}
{n \choose r}H_r(p)=S_r(p)=\sigma_r(\lambda_1(p),...,\lambda_n(p)),
\end{eqnarray}
where $\sigma_r$ is the $r$-th symmetric function on $\mathbb R^n$,
\begin{eqnarray}
\sigma_r(x_1,\ldots,x_n)=\sum_{i_1<\ldots<i_r}x_{i_1}\ldots x_{i_r}.
\end{eqnarray}
For future use, we observe that
\begin{eqnarray}\label{hr9}
\sigma_r(x)=x_i\sigma_{r-1}(\widehat{x_i})+\sigma_r(\widehat{x_i}),\;\;\;\;i,r=1,\ldots,n,\;\;x\in\mathbb R^n,
\end{eqnarray}
where $\sigma_0=1$ and, for instance, $\sigma_{r-1}(\widehat{x_i})=\sigma_{r-1}(x_1,\ldots,x_{i-1},x_{i+1},\ldots,x_n)$.

The squared norm $|A|^2$ of the shape operator $A$ is defined as the trace of $A^2$. Taking an orthonormal basis that diagonalizes $A$, it is easy to see that
\begin{eqnarray}\label{tb116}
|A|^2=\sum_{i}\lambda_i^2.
\end{eqnarray}
Hence, from (\ref{Hr}), (\ref{hr2}) and (\ref{tb116}), we obtain
\begin{eqnarray}\label{tb120}
n^2H^2=\Big(\sum_{i=1}^n\lambda_i\Big)^2=\sum_{i=1}^n\lambda_i^2+\sum_{i\neq j}\lambda_i\lambda_j=|A|^2+n(n-1)(R-c).
\end{eqnarray}

We will finish this section recalling the definition of the Newton's tensors $P_r,\;0\leq r\leq n$, associated to the shape operator $A$ of an isometric immersion. These tensors are defined inductively by
\begin{eqnarray}\label{hr16}
P_0&=&I,\nonumber\\P_r&=&S_rI-AP_{r-1},\;\;\;1\leq r\leq n.
\end{eqnarray}

\section{Main ingredients}

In this section we will present the basic tools for the proofs of our results. We begin stating the following proposition \cite[p. 668]{li}, which provides a formula for the Laplacian of $|A|^2$ for a hypersurface in a space of constant curvature.

\begin{prop}\label{prop_simons}
Let $f:M^n\to N_c^{n+1}$ be an isometric immersion, $p\in M^n$. Denote by $\lambda_1,\ldots,\lambda_n$ the principal curvatures of $M^n$ at $p$ and let ${e_1,\ldots,e_n}$ be an orthonormal basis of $T_p M$ such that $Ae_i=\lambda_ie_i$, $i=1,\ldots,n$. Then
\begin{eqnarray}\label{formula_simons}
\frac{1}{2}\triangle|A|^2=|\nabla A|^2+n\sum_i \lambda_i\textnormal{Hess}\,H(e_i,e_i)+\sum_{i<j}(\lambda_i-\lambda_j)^2K_{ij},
\end{eqnarray}
where $|\nabla A|$ is the norm of the covariant derivative $\nabla A$ of $A$, $\textnormal{Hess}\,H$ is the Hessian of $H$ and $K_{ij}$ is the sectional curvature of $M^n$ in the plane spanned by $\{e_i,e_j\}$.
\end{prop}

Using (\ref{formula_simons}) and the Gauss equation, one easily obtains
\begin{eqnarray}\label{formula_simons2}
\frac{1}{2}\triangle|A|^2=|\nabla A|^2+n\sum_i \lambda_iH_{ii}+nc(|A|^2-nH^2)+nH\text{tr}(A^3)-|A|^4,
\end{eqnarray}
where $\text{tr}(A^3)$ is the trace of $A^3$ and $H_{ii}=\text{Hess}\,H(e_i,e_i)$.

\vskip10pt

In this paper we deal with hypersurfaces with constant mean curvature $H$ in Euclidean spaces. When working with hypersurfaces of constant mean constant, it is convenient to replace the shape operator $A$ by the tensor field $\phi=A-HI$. It is easy to see that $\phi$ is symmetric and that its eigenvalues are $\mu_i=\lambda_i-H$, $i=1,\ldots,n$, where $\lambda_1,\ldots,\lambda_n$ are the eigenvalues of $A$. It is easy to see that
\begin{eqnarray}
\text{tr}(\phi)&=&0,\label{phi1}\\
|\phi|^2&=&|A|^2-nH^2,\label{phi2}\\
\text{tr}(A^3)&=&\text{tr}(\phi^3)+3H|\phi|^2+nH^3,\label{phi3}
\end{eqnarray}
where $|\phi|^2=\text{tr}(\phi^2)$. From (\ref{phi2}) one has that $|\phi|\equiv0$
if and only if the immersion is totally umbilical.

In our proofs we will also need of the following result of Okumura \cite{Ok}:
\begin{lem}[Okumura]\label{okumura}
Let $\mu_1,\mu_2,\ldots,\mu_n$, $n\geq 3$, be real numbers such that $\sum_i\mu_i=0$ and $\sum_i\mu_i^2=\beta^2$. Then
\begin{eqnarray}\label{des_okumura}
-\frac{n-2}{\sqrt{n(n-1)}}\beta^3\leq\sum_i\mu_i^3\leq\frac{n-2}{\sqrt{n(n-1)}}\beta^3,
\end{eqnarray}
and equality holds in (\ref{des_okumura}) if and only if at least $(n-1)$ of the $\mu_i$ are equal.
\end{lem}
Let $f:M^n\to N_c^{n+1}$ be an isometric immersion with constant mean curvature $H$. From (\ref{formula_simons2}), (\ref{phi2}) and (\ref{phi3}), one obtains
\begin{eqnarray*}
\frac{1}{2}\triangle|A|^2&=&|\nabla A|^2+nc(|A|^2-nH^2)+nH\text{tr}(A^3)-|A|^4\\
&=&|\nabla A|^2+nc|\phi|^2+
nH\big[\text{tr}(\phi^3)+3H|\phi|^2+nH^3\big]\\
&&-\big(|\phi|^2+nH^2\big)^2\\
&=&|\nabla A|^2+nH\text{tr}(\phi^3)+|\phi|^2[nc+
nH^2-|\phi|^2].
\end{eqnarray*}
Hence, by Lema \ref{okumura},
\begin{eqnarray}\label{simons phi}
\frac{1}{2}\triangle|A|^2\geq |\nabla A|^2+|\phi|^2\Big[nc+
nH^2-\frac{n(n-2)}{\sqrt{n(n-1)}}|H||\phi|-|\phi|^2\Big].
\end{eqnarray}

Another important tool for the proofs of our results is the following well known result of Smyth and Xavier \cite{Smyth_Xa}.

\vskip10pt

\noindent{\bf The principal curvature theorem.}
{\it Let $M^n$ be a complete orientable Riemannian manifold , and $f:M^n\to\mathbb R^{n+1}$ be a isometric immersion such that $f(M^n)$ is not a hyperplane. Let $\Lambda\subset\mathbb R$ be the set of nonzero values assumed by the principal curvatures of shape operator $A$ and let $\Lambda^{\pm}=\Lambda\cap\mathbb R^{\pm}$.
\begin{enumerate}
\item [(i)] If $\Lambda^{+}\neq\emptyset$ and $\Lambda^{-}\neq\emptyset$, then $\emph{inf}\Lambda^{+}=0=\emph{sup}\Lambda^{-}$.
\item [(ii)] If $\Lambda^{+}$ or $\Lambda^{-}$ is empty, then $\overline{\Lambda}$ is connected.
\end{enumerate}}

We will finish this section stating a result that classifies the complete isoparametric hypersurfaces in an Euclidean space. Recall that a hypersurface of the complete simply-connected $(n+1)$-dimensional space form $\mathbb Q^{n+1}_c$ of constant sectional curvature $c$ is called isoparametric when its principal curvatures are constant. The classification of such hypersurfaces in the case $c=0$ was established by Levi-Civita \cite{leviCivita} for $n=2$ and by Segre \cite{Segre} for higher dimensions. By their works, one has the following result:

\begin{thm}\label{isopar_euclid}
If $M^n$ is a complete connected Riemannian manifold and $f:M^n\to\mathbb R^{n+1}$ an isoparametric isometric immersion, then $f(M^n)$ is a generalized cylinder $\mathbb R^{n-k}\times\mathbb S^{k}(r)$, for some $k=0,\ldots,n$ and some $r>0$.
\end{thm}

\section{Proof of Theorem \ref{teorema1}}\label{prova_Teo1}
Choose the orientation so that $H>0$. Denoting by $\lambda_1\leq\lambda_2\leq\lambda_3\leq\lambda_4$ the principal curvatures of $M^4$, one has, at each point of $M^4$,
\begin{eqnarray}\label{q1}
\lambda_4\geq H>0.
\end{eqnarray}
Since, by hypothesis, $H$ and $R$ are constant, it follows from (\ref{tb120}) that $|A|^2$ is constant. Hence, by  (\ref{tb117}) and Proposition \ref{prop_simons},
\begin{eqnarray}\label{q2}
0=\frac{1}{2}\Delta|A|^2&=&|\nabla A|^2+\sum_{i<j}(\lambda_i-\lambda_j)^2\lambda_i \lambda_j\nonumber\\
&\geq&\sum_{i<j}(\lambda_i-\lambda_j)^2\lambda_i \lambda_j.
\end{eqnarray}
\noindent{\bf Claim.} $\inf|\lambda_3|>0$.\vspace{0.25cm}

Suppose, by contradiction, that there is a sequence $(p_k)$ in $M$ such that $\lambda_3(p_k)\to0$. Since the principal curvature functions are bounded (because $|A|$ is constant), passing to a subsequence if necessary, we can assume that $\lambda_i(p_k)\to\overline{\lambda}_i$, $\forall i$. By (\ref{q1}),
\begin{eqnarray}\label{q3}
\overline{\lambda}_1\leq\overline{\lambda}_2\leq\overline{\lambda}_3=0<H\leq\overline{\lambda}_4.
\end{eqnarray}
Since $H$ and $R$ are constant, one has
\begin{eqnarray*}
4H&=&\overline{\lambda}_1+\overline{\lambda}_2+\overline{\lambda}_4\\
6R&=&\overline{\lambda}_1\overline{\lambda}_2+\overline{\lambda}_1\overline{\lambda}_4+\overline{\lambda}_2\overline{\lambda}_4.
\end{eqnarray*}
Multiplying the first of the equalities above by $\overline{\lambda}_2$ and using the second, we obtain, with the aid of (\ref{q3}),
$$4H\overline{\lambda}_2=\overline{\lambda}_1\overline{\lambda}_2+\overline{\lambda}_2^2
+\overline{\lambda}_2\overline{\lambda}_4=6R-\overline{\lambda}_1\overline{\lambda}_4+\overline{\lambda}_2^2\geq6R>0,$$
whence one obtains $\overline{\lambda}_2>0$. Since this contradicts (\ref{q3}), the claim is proved.\\

It follows from the claim that $\sup\lambda_3<0$ or $\inf\lambda_3>0$. If we had $\sup\lambda_3<0$, then $\sup\Lambda^{-}=\sup\lambda_3<0$, contradicting the principal curvature theorem. Hence,
\begin{eqnarray}\label{q4}
\inf\lambda_3>0.
\end{eqnarray}

Assume first that there is a sequence $(p_k)$ in $M$ such that $\lambda_2(p_k)\to0$. Passing to a subsequence, we can assume that $\lambda_i(p_k)\to\overline{\lambda}_i$, $\forall i$. By (\ref{q4}),
\begin{eqnarray}\label{q5}
\overline{\lambda}_1\leq\overline{\lambda}_2=0<\overline{\lambda}_3\leq\overline{\lambda}_4.
\end{eqnarray}
Since $H$ and $R$ are constant, one has
\begin{eqnarray}\label{q6}
4H&=&\overline{\lambda}_1+\overline{\lambda}_3+\overline{\lambda}_4,\label{nz413}\\
6R&=&\overline{\lambda}_1 \overline{\lambda}_3+\overline{\lambda}_1 \overline{\lambda}_4+\overline{\lambda}_3\overline{\lambda}_4.\nonumber
\end{eqnarray}
Multiplying the first of the equalities above by $\overline{\lambda}_4$ and using the second, we obtain
\begin{eqnarray*}
4H\overline{\lambda}_4=\overline{\lambda}_1\overline{\lambda}_4+\overline{\lambda}_3\overline{\lambda}_4+\overline{\lambda}_4^2
=6R-\overline{\lambda}_1\overline{\lambda}_3+\overline{\lambda}_4^2.
\end{eqnarray*}
It follows from (\ref{q5}), the above equality and the hypothesis $R\geq\frac{2}{3}H^2$ that
\begin{eqnarray*}
0\geq\overline{\lambda}_1\overline{\lambda}_3=\overline{\lambda}_4^2-4H\overline{\lambda}_4+6R
&=&(\overline{\lambda}_4-2H)^2-4H^2+6R\\
&=&(\overline{\lambda}_4-2H)^2+6\Big(R-\frac{2}{3}H^2\Big)\\
&\geq&0,
\end{eqnarray*}
and so
\begin{eqnarray}\label{q7}
\overline{\lambda}_1=0,\;\;\;\overline{\lambda}_4=2H\;\;\;\text{e}\;\;\;R=\frac{2}{3}H^2.
\end{eqnarray}
Hence, by (\ref{q5}), (\ref{q6}) and (\ref{q7}),
$$\overline{\lambda}_1=\overline{\lambda}_2=0\;\;\;\;\text{e}\;\;\;\overline{\lambda}_3=\overline{\lambda}_4=2H.$$

Suppose now that $\inf|\lambda_2|>0$. By (\ref{q4}) and the principal curvature theorem, one has
\begin{eqnarray}\label{q8}
\inf\lambda_2>0.
\end{eqnarray}
We have two possibilities:
\begin{enumerate}
\item [1)] There exists some $p\in M$ such that $\lambda_1(p)=0$.
\item [2)] $\lambda_1(x)\neq 0$, $\forall x\in M$.
\end{enumerate}

Assuming 1), one has, by (\ref{q2}) and (\ref{q8}),
\begin{eqnarray}\label{q9}
\lambda_1(p)=0,\;\;\lambda_2(p)=\lambda_3(p)=\lambda_4(p)=\frac{4}{3}H.
\end{eqnarray}
Hence, by (\ref{tb116}), (\ref{tb120}) and (\ref{phi2}),
$$|A|^2=\frac{16}{3}H^2,\;\;\;R=\frac{8}{9}H^2\;\;\;\;\text{e}\;\;\;|\phi|^2=\frac{4}{3}H^2.$$
Substituting in (\ref{simons phi}), we arrive at
\begin{eqnarray*}
0=\frac{1}{2}\triangle|A|^2&\geq& |\nabla A|^2+|\phi|^2\Big(
4H^2-\frac{4}{\sqrt{3}}H|\phi|-|\phi|^2\Big)\\
&=&|\nabla A|^2+|\phi|^2\Big(
4H^2-\frac{8}{3}H^2-\frac{4}{3}H^2\Big)\\
&=&|\nabla A|^2.
\end{eqnarray*}
Hence, $\nabla A\equiv 0$ and therefore $f$ is isoparametric \cite[p. 254]{RYAN}. Now it follows from (\ref{q9}) and Theorem \ref{isopar_euclid} that
$$f(M^4)=\mathbb R\times\mathbb S^3\Big(\frac{3}{4H}\Big).$$

Assuming 2), one has $\lambda_1>0$ along $M$. Indeed, if we had $\lambda_1(q)<0$ at some point $q\in M$, then by continuity $\lambda_1$ would be negative along $M$, and from (\ref{q8}) one would obtain $\inf\Lambda^{+}>0$, contradicting the principal curvature theorem.
From $\lambda_1>0$ and (\ref{q2}), we obtain
\begin{eqnarray}\label{q10}
\lambda_1\equiv\lambda_2\equiv\lambda_3\equiv\lambda_4\equiv H,
\end{eqnarray}
that is, $f$ is totally umbilical. It now follows from the classification of umbilical hypersurfaces in an Euclidean space that $f(M)$ is contained in a hypersphere or in a hyperplane of $\mathbb R^5$. Since $M$ is complete and $H>0$, $f(M)$ is a hypersphere of $\mathbb R^5$. Hence, by (\ref{q10}), $R=H^2$ and
$$f(M^4)=\mathbb S^4\Big(\frac{1}{H}\Big).\qed$$

\section{Proof of Theorem \ref{teorema2}}\label{prova_Teo2}

Choose the orientation so that $H>0$. Denoting by $\lambda_1\leq\lambda_2\leq\lambda_3\leq\lambda_4\leq\lambda_5$ the principal curvatures of $M^5$, one has, at each point of $M^5$,
\begin{eqnarray}\label{bn1}
\lambda_5\geq H>0.
\end{eqnarray}
Since, by hypothesis, $H$ and $R$ are constant, it follows from (\ref{tb120}) that $|A|^2$ is constant. Hence, by (\ref{tb117}) and Proposition \ref{prop_simons},
\begin{eqnarray}\label{bn2}
0=\frac{1}{2}\triangle|A|^2&=&|\nabla A|^2+\sum_{i<j}(\lambda_i-\lambda_j)^2\lambda_i \lambda_j\nonumber\\
&\geq&\sum_{i<j}(\lambda_i-\lambda_j)^2\lambda_i \lambda_j.
\end{eqnarray}
By (\ref{hr2}), (\ref{hr3}) and (\ref{hr9}), one has
\begin{eqnarray*}
10R=S_2=\lambda_i\sigma_1(\widehat{\lambda_i})+\sigma_2(\widehat{\lambda_i})
=\lambda_i(5H-\lambda_i)+\sigma_2(\widehat{\lambda_i}),\;\;\forall i.
\end{eqnarray*}
From the above equality and the hypothesis $R\geq\frac{5}{8}H^2$, one obtains
\begin{eqnarray}\label{bn3}
\sigma_2(\widehat{\lambda_i})&=&\lambda_i^2-5H\lambda_i+10R\nonumber\\
&=&\Big(\lambda_i-\frac{5H}{2}\Big)^2-\frac{25}{4}H^2+10R\nonumber\\
&=&\Big(\lambda_i-\frac{5H}{2}\Big)^2+10\Big(R-\frac{5}{8}H^2\Big)\nonumber\\
&\geq&0,\;\;\forall i.
\end{eqnarray}
\noindent{\bf Claim.} $\inf|\lambda_4|>0$.\vspace{0.25cm}

Suppose, by contradiction, that there is a sequence $(p_k)$ in $M$ such that $\lambda_4(p_k)\to0$. Since the principal curvature functions are bounded (because $|A|$ is constant), passing to a subsequence if necessary, we can assume that $\lambda_i(p_k)\to\overline{\lambda}_i$, $\forall i$. By (\ref{bn1}),
\begin{eqnarray}\label{bn3a}
\overline{\lambda}_1\leq\overline{\lambda}_2\leq\overline{\lambda}_3
\leq\overline{\lambda}_4=0<H\leq\overline{\lambda}_5.
\end{eqnarray}
Since $\sigma_4(\overline{\lambda}_1,\overline{\lambda}_2,\overline{\lambda}_3,\overline{\lambda}_4,\overline{\lambda}_5)\geq 0$ (because, by hypothesis, $H_4\geq0$), one has $\overline{\lambda}_3=0$ and therefore,
\begin{eqnarray*}
0<10R=\overline{\lambda}_1\overline{\lambda}_2+\overline{\lambda}_1\overline{\lambda}_5+\overline{\lambda}_2\overline{\lambda}_5
\leq\overline{\lambda}_1\overline{\lambda}_2,
\end{eqnarray*}
where in the last inequality we used (\ref{bn3a}). Hence,
\begin{eqnarray}\label{bn4}
\overline{\lambda}_1\leq\overline{\lambda}_2<\overline{\lambda}_3=\overline{\lambda}_4=0<H\leq\overline{\lambda}_5,
\end{eqnarray}
which implies
\begin{eqnarray}\label{bn5}
\lim_{k\to\infty}S_3(p_k)=\sigma_3(\overline{\lambda}_1,\overline{\lambda}_2,\overline{\lambda}_3,\overline{\lambda}_4,\overline{\lambda}_5)
=\overline{\lambda}_1\overline{\lambda}_2\overline{\lambda}_5>0.
\end{eqnarray}
On the other hand, by (\ref{hr9}), (\ref{bn3}) and (\ref{bn4}), one has
\begin{eqnarray*}
\lim_{k\to\infty}S_3(p_k)&=&\overline{\lambda}_1\sigma_2(\overline{\lambda}_2,\overline{\lambda}_3,\overline{\lambda}_4,\overline{\lambda}_5)
+\sigma_3(\overline{\lambda}_2,\overline{\lambda}_3,\overline{\lambda}_4,\overline{\lambda}_5)\nonumber\\
&\leq&\sigma_3(\overline{\lambda}_2,\overline{\lambda}_3,\overline{\lambda}_4,\overline{\lambda}_5)\nonumber\\
&=&0,
\end{eqnarray*}
contradicting (\ref{bn5}). This proves the claim.

From the claim and the principal curvature theorem we obtain
\begin{eqnarray}\label{bn6}
\inf\lambda_4>0.
\end{eqnarray}

Assume first that there is a sequence $(p_k)$ in $M$ such that $\lambda_3(p_k)\to0$. Passing to a subsequence, we can assume that $\lambda_i(p_k)\to\overline{\lambda}_i$, $\forall i$. By (\ref{bn6}),
\begin{eqnarray}\label{bn6b}
\overline{\lambda}_1\leq\overline{\lambda}_2\leq\overline{\lambda}_3=0<\overline{\lambda}_4\leq\overline{\lambda}_5,
\end{eqnarray}
and by (\ref{bn2}),
\begin{eqnarray}\label{bn6a}
0\geq\sum_{i<j}(\overline{\lambda}_i-\overline{\lambda}_j)^2\overline{\lambda}_i\overline{\lambda}_j.
\end{eqnarray}
From (\ref{hr9}), (\ref{bn3}) and (\ref{bn6b}), one obtains
\begin{eqnarray}\label{bn7}
\lim_{k\to\infty}S_3(p_k)&=&\overline{\lambda}_5\sigma_2(\overline{\lambda}_1,\overline{\lambda}_2,\overline{\lambda}_3,\overline{\lambda}_4)
+\sigma_3(\overline{\lambda}_1,\overline{\lambda}_2,\overline{\lambda}_3,\overline{\lambda}_4)\nonumber\\
&\geq&\sigma_3(\overline{\lambda}_1,\overline{\lambda}_2,\overline{\lambda}_3,\overline{\lambda}_4)
=\overline{\lambda}_1\overline{\lambda}_2\overline{\lambda}_4\geq0.
\end{eqnarray}
On the other hand,
\begin{eqnarray}\label{bn8}
\lim_{k\to\infty}S_3(p_k)&=&\overline{\lambda}_2\sigma_2(\overline{\lambda}_1,\overline{\lambda}_3,\overline{\lambda}_4,\overline{\lambda}_5)+
\sigma_3(\overline{\lambda}_1,\overline{\lambda}_3,\overline{\lambda}_4,\overline{\lambda}_5)\nonumber\\
&\leq&\sigma_3(\overline{\lambda}_1,\overline{\lambda}_3,\overline{\lambda}_4,\overline{\lambda}_5)
=\overline{\lambda}_1\overline{\lambda}_4\overline{\lambda}_5\leq 0.
\end{eqnarray}
It follows from (\ref{bn6b}), (\ref{bn7}) and (\ref{bn8}) that $\overline{\lambda}_1=0$. Hence, by (\ref{bn6b}) and (\ref{bn6a}),
\begin{eqnarray*}
\overline{\lambda}_1=\overline{\lambda}_2=\overline{\lambda}_3=0\;\;\;\text{e}\;\;\overline{\lambda}_4=\overline{\lambda}_5=\frac{5H}{2},
\end{eqnarray*}
and so
\begin{eqnarray*}
R=\frac{5}{8}H^2.
\end{eqnarray*}

Suppose now that $\inf|\lambda_3|>0$. By (\ref{bn6}) and the principal curvature theorem, one has
\begin{eqnarray}\label{bn10}
\inf\lambda_3>0.
\end{eqnarray}

Assume that there is a sequence $(p_k)$ in $M$ such that $\lambda_2(p_k)\to0$. Passing to a subsequence, we can assume that $\lambda_i(p_k)\to \overline{\lambda}_i$, $\forall i$. By (\ref{bn10}),
\begin{eqnarray}\label{bn11}
\overline{\lambda}_1\leq\overline{\lambda}_2=0<\overline{\lambda}_3\leq\overline{\lambda}_4\leq\overline{\lambda}_5.
\end{eqnarray}
From (\ref{bn2}), we obtain
\begin{eqnarray}\label{bn12}
0\geq\sum_{i<j}(\overline{\lambda}_i-\overline{\lambda}_j)^2\overline{\lambda}_i\overline{\lambda}_j.
\end{eqnarray}
Since $\sigma_4(\overline{\lambda}_1,\overline{\lambda}_2,\overline{\lambda}_3,\overline{\lambda}_4,\overline{\lambda}_5)\geq0$, one has $\overline{\lambda}_1=0$. It follows from (\ref{bn11}) and (\ref{bn12}) that
\begin{eqnarray*}
\overline{\lambda}_1=\overline{\lambda}_2=0\;\;\;\text{e}\;\;\overline{\lambda}_3=\overline{\lambda}_4=\overline{\lambda}_5=\frac{5}{3}H,
\end{eqnarray*}
and so
$$R=\frac{5}{6}H^2.$$

Assume now that $\inf|\lambda_2|>0$. By (\ref{bn10}) and the principal curvature theorem, one has
\begin{eqnarray}\label{bn13}
\inf\lambda_2>0.
\end{eqnarray}
We have two possibilities:
\begin{enumerate}
\item [1)] There exists some $p\in M$ such that $\lambda_1(p)=0$.
\item [2)] $\lambda_1(x)\neq0$, $\forall x\in M$.
\end{enumerate}

Assuming 1), one has, by (\ref{bn2}) and (\ref{bn13}),
\begin{eqnarray}\label{bn14}
\lambda_1(p)=0,\;\;\lambda_2(p)=\lambda_3(p)=\lambda_4(p)=\lambda_5(p)=\frac{5}{4}H.
\end{eqnarray}
Hence, by (\ref{tb116}), (\ref{tb120}) and (\ref{phi2}),
\begin{eqnarray*}
|A|^2=\frac{25}{4}H^2,\;R=\frac{15}{16}H^2\;\;\;\text{e}\;\;\;|\phi|^2=\frac{5}{4}H^2.
\end{eqnarray*}
Substituting in (\ref{simons phi}), we arrive at
\begin{eqnarray*}
0=\frac{1}{2}\triangle|A|^2&\geq& |\nabla A|^2+|\phi|^2\Big(
5H^2-\frac{15}{2\sqrt{5}}H|\phi|-|\phi|^2\Big)\\
&=&|\nabla A|^2+|\phi|^2\Big(
5H^2-\frac{15}{4}H^2-\frac{5}{4}H^2\Big)\\
&=&|\nabla A|^2.
\end{eqnarray*}
Hence, $\nabla A\equiv0$ and therefore $f$ is isoparametric \cite[p. 254]{RYAN}. It now follows from (\ref{bn14}) and Theorem \ref{isopar_euclid} that
$$f(M^5)=\mathbb R\times\mathbb S^4\Big(\frac{4}{5H}\Big).$$

Assuming 2), one has $\lambda_1>0$ along $M$. Indeed, if we had $\lambda_1(q)<0$ at some point $q\in M$, then by continuity $\lambda_1$ would be negative along $M$, and from (\ref{bn13}) one would obtain $\inf\Lambda^{+}>0$, contradicting the principal curvature theorem.
From $\lambda_1>0$ and (\ref{bn2}), we obtain
\begin{eqnarray}\label{bn15}
\lambda_1\equiv\lambda_2\equiv\lambda_3\equiv\lambda_4\equiv\lambda_5\equiv H,
\end{eqnarray}
that is, $f$ is totally umbilical. It now follows from the classification of umbilical hypersurfaces in an Euclidean space that $f(M)$ is contained in a hypersphere or in a hyperplane of $\mathbb R^6$. Since $M$ is complete and $H>0$, $f(M)$ is  a hypersphere of $\mathbb R^6$. Hence, by (\ref{bn15}), $R=H^2$ and $$f(M^5)=\mathbb S^5\Big(\frac{1}{H}\Big).\qed$$

\section{Hypersurfaces of arbitrary dimension}
In the proof of Theorem \ref{teorema3} we will use the following lemma:

\begin{lem}\label{trA3_lem}
Let $M^n$ be a orientable Riemannian manifold $M^n$, with dimension $n\geq3$, and let $f:M^n\to N_c^{n+1}$ be an isometric immersion. Denote by $A$ the shape operator of the immersion with respect to a global unit normal vector $\xi$. Then,
\begin{eqnarray}\label{trA3}
\emph{tr} A^3=\frac{nH}{2}\big(3|A|^2-n^2H^2\big)+3S_3.
\end{eqnarray}
\end{lem}
\noindent{\bf Proof.} It is well known (see, for example, \cite[Lema 2.1.]{BarCol}) that
\begin{eqnarray*}
\text{tr}(AP_r)=(r+1)S_{r+1},\;1\leq r\leq n-1,
\end{eqnarray*}
where $P_r$ is the $r$-th Newton's tensor (cf. Section \ref{preliminar}). Making $r=2\;$  in the above equality, one has, by (\ref{hr16}),
\begin{eqnarray}\label{rt1}
3S_3=\text{tr}(AP_2)=\text{tr}(S_2A-S_1A^2+A^3)=S_2nH-S_1|A|^2+\text{tr}A^3.
\end{eqnarray}
By (\ref{hr2}), (\ref{hr3}) and (\ref{tb120}), we have
\begin{eqnarray}\label{rt2}
n^2H^2=|A|^2+2S_2.
\end{eqnarray}
From (\ref{rt1}) and (\ref{rt2}), we obtain
\begin{eqnarray*}
\text{tr}A^3&=&-nH\Big(\frac{n^2H^2}{2}-\frac{|A|^2}{2}\Big)+nH|A|^2+3S_3\\
&=&\frac{nH}{2}\big(3|A|^2-n^2H^2\big)+3S_3.\qed
\end{eqnarray*}

\noindent{\bf Proof of Theorem \ref{teorema3}.}
Since, by hypothesis, $H$ and $R$ are constant, one has by (\ref{tb120}) that $|A|^2$ is also constant. Hence, by (\ref{formula_simons2}) and Lemma \ref{trA3_lem} one has, since $HH_3\geq 0$ and $0\leq R\leq \frac{nH^2}{2(n-1)}$ by hypothesis,
\begin{eqnarray}\label{formula_BB}
0=\frac{1}{2}\triangle|A|^2&=&|\nabla A|^2+\frac{n^2H^2}{2}(3|A|^2-n^2H^2)+3nHS_3-|A|^4\nonumber\\
&=&|\nabla A|^2-\Big[|A|^4-\frac{3n^2H^2}{2}|A|^2+\frac{(n^2H^2)^2}{2}\Big]+3nHS_3\nonumber\\
&=&|\nabla A|^2-\Big(|A|^2-\frac{n^2H^2}{2}\Big)\Big(|A|^2-n^2H^2\Big)+3nHS_3\nonumber\\
&=&|\nabla A|^2+n^2(n-1)^2R\Big(\frac{nH^2}{2(n-1)}-R\Big)+3nHS_3\nonumber\\
&\geq&|\nabla A|^2.
\end{eqnarray}
Therefore, $\nabla A\equiv0$ and $R=0$ or $R=\frac{nH^2}{2(n-1)}$. From $\nabla A\equiv0$ one concludes that $f$ is isoparametric \cite[p. 254]{RYAN}. Since $n\geq 3$ and $H\neq 0$, it follows from Theorem \ref{isopar_euclid} that
$$f(M^n)=\mathbb R^{k}\times\mathbb S^{n-k}(r),$$
for some $r>0$ and some $1\leq k\leq n-1$. Using (\ref{tb120}), it can be easily verified that $R=0$ occurs precisely when $k=n-1$, and $R=\frac{nH^2}{2(n-1)}$
when $k=n-2$. Hence, $f(M^n)=\mathbb R^{n-1}\times\mathbb S^1\Big(\frac{1}{n|H|}\Big)$ when $R=0$, and $f(M^n)=\mathbb R^{n-2}\times\mathbb S^2\Big(\frac{2}{n|H|}\Big)$ when $R=\frac{nH^2}{2(n-1)}$.
\qed

\noindent Roberto Alonso N\'u\~nez\\Rua Dr. Paulo Alves 110, Bl C, Apto. 402\\24210-445 Niter\'oi, RJ, Brazil
\\\texttt{\detokenize{alonso_nunez@id.uff.br}}


\begin{thebibliography}{s2}

\bibitem{BarCol} {\sc J.L.M. Barbosa and A.G. Colares}, {\em Stability of hypersurfaces with constant r-mean curvature},
Ann. Global Anal. Geom. {\bf 15} (1997), 277-297.

\bibitem{ChengWan} {\sc Q.M. Cheng and Q.R. Wan}, {\em Complete hypersurfaces of $\mathbb R^4$ with constant mean curvature},
Monatsh. Math. {\bf 118} (1994), 171-204.

\bibitem{ChengYauMean} {\sc S. Y. Cheng and S. T. Yau}, {\em Differential Equations on Riemannian Manifolds and
their Geometric Applications},
Comm. Pure Appl. Math., {\bf 28} (1975) 333-354.

\bibitem{ChengYauScal} {\sc S.Y. Cheng and S.T. Yau}, {\em Hypersurfaces with constant scalar curvature}, Math.
Ann., {\bf 225} (1977), 195-204.

\bibitem{HART} {\sc P. Hartman}, {\em On complete hypersurfaces of nonnegative sectional curvatures and constant mth mean curvature}, Trans. Amer. Math. Soc., {\bf{245}} (1978) 363-374.

\bibitem{HUANG} {\sc X.G. Huang}, {\it Complete Hypersurfaces with constant scalar curvature and constant mean curvature in $\mathbb R^4$},
Chin. Ann. Math. Ser. B, Vol. {\bf{6B(2)}} (1985) 177-184.

\bibitem{Klot} {\sc T. Klotz and R. Osserman}, {\em Complete surfaces in $E^3$ with constant mean curvature}, Comm. Math. Helv. {\bf{41}} (1966/67) 313-318.

\bibitem{leviCivita} {\sc T. Levi-Civita}, {\em Famiglie di superficie isoparametriche nell'ordinario spazio euclideo}, Atti Accad. Naz. Lincei Rend. Cl. Sci. Fis. Mat. Natur. (6) {\bf 26} (1937) 355-362.

\bibitem{li} {\sc H. Li}, {\em Hypersurfaces with constant scalar curvature in space forms}, Math. Ann. {\bf 305} (1996) 665-672.

\bibitem{Nom_Smyth} {\sc K. Nomizu and B. Smyth}, {\it A formula of Simon's type and hypersurfaces with constant mean curvature}, J. Differential Geom., {\bf 3} (1969) 367-377.

\bibitem{Ok} {\sc M. Okumura}, {\em Hypersurfaces and pinching problem on the second fundamental tensor},
Amer. J. Math., {\bf 96} (1974) 207-213.

\bibitem{RYAN} {\sc P.J. Ryan}, {\em Hypersurfaces with parallel Ricci tensor}, Osaka J. Math. {\bf 8} (1971) 251-259.

\bibitem{Segre} {\sc B. Segre}, {\em Famiglie di ipersuperficie isoparametriche negli spazi euclidei ad un qualunque numero di dimensioni}, Atti Accad. Naz. Lincei Rend. Cl. Sci. Fis. Mat. Natur. (6) {\bf 27} (1938) 203-207.

\bibitem{Smyth_Xa} {\sc B. Smyth and F. Xavier}, {\em Efimov's theorem in dimension greater than two},
§Invent. Math., {\bf 90} (1987) 443-450.

\end{thebibliography}
\end{document}